\documentclass[leqno,11pt]{amsart}

\usepackage{amssymb}
\usepackage{color}
\usepackage{latexsym}
\usepackage{amscd}
\usepackage{enumerate} 
\theoremstyle{plain}
\newcounter{theo}
\newtheorem{thm}[theo]{Theorem}
\newtheorem{prop}[theo]{Proposition}
\newtheorem{lem}[theo]{Lemma}
\newtheorem{exa}[theo]{Example}
\theoremstyle{definition}
\newtheorem{rem}[theo]{Remark}
\newtheorem{df}[theo]{Definition}
\newenvironment{pf}{{\noindent\bf Proof. }}{\hfill $\square$\medskip}
\title[Nondegeneracy conditions]{On  nondegeneracy conditions  for the Levi map in higher codimension: a Survey}
\author{L\'ea Blanc-Centi and Francine Meylan}

\newcommand{\C}{\mathbb{C}}
\newcommand{\R}{\mathbb{R}}
\newcommand{\Aut}{\mathrm{Aut}}
\newcommand{\Imm}{\mathrm{Im}\,}
\newcommand{\Ree}{\mathrm{Re}\,}
\newcommand{\transp}{\,^t}
\newcommand{\aaa}{\mathfrak{(a)}}
\newcommand{\bbb}{\mathfrak{(b)}}

\begin{document}

\maketitle

\begin{abstract}
We compare various definitions of nondegeneracy  of the Levi map for real submanifolds of higher codimension in $\C^N$ and discuss the generalization to higher codimension  of the   $2$-jet determination for biholomorphisms  in the hypersurface case proved by Chern and Moser in \cite{CM}.
\end{abstract}

\section{Introduction}
Let $M$ be a real submanifold of $\C^N$, $p\in M$, and $\Aut(M,p)$ be the stability group of $M$ at point $p$, that is the set of (germs of) biholomorphisms $F$ fixing $p$ and such that $F(M)\subset M$. 

For a real hypersurface in complex dimension 2, H. Poincar\'e initiated the study of the stability group  by looking at Taylor series expansion: the condition $F(M)\subset M$ means that $\rho(F(z,w))_{|M}=0$, where $\rho$ is a defining function of $M$, and this equation gives some constraints on the Taylor series coefficients of $F$. The process was carried out much later in a significant manner by J.K. Moser for Levi non-degenerate hypersurfaces, to obtain the following $2$-jet determination statement:

\begin{thm} \cite{CM}
Let $M$ be a real-analytic hypersurface through a point $p$ in $\C^N$ with non-degenerate Levi form at $p$. Let $F$, $G$ be two germs of biholomorphic maps preserving $M$. Then, if $F$ and $G$ have the same 2-jets at $p$, they coincide.
\end{thm}
\noindent The proof relies on the fact that the elements of $\Aut(M,p)$ have to satisfy an overdeterminated inhomogeneous system of partial differential equations, which is a very restrictive condition. But the result becomes false without any hypothesis on the Levi form, as one can see by considering the hyperplane $\mathrm{Im}\,w=0$ in $\C_z^{N}\times\C_w$, whose stability group at 0 is infinite dimensional.

\bigskip

Here, we are interested in $2$- jet determination problems in higher codimension  in the context of "nondegeneracy" of the Levi map.  For $M$ being a submanifold of codimension $d$ in $\C^N$, the first step is to choose a convenient generalization of the Levi nondegeneracy condition. Various definitions appeared in the literature. Our aim is to survey these definitions, to compare them by giving many examples, and to explain why the definition introduced by Beloshapka in \cite{Be2} seems the most relevant for us in the approach by model submanifolds. More precisely, we present a detailed explanation of the first part of his original  proof in \cite{Be2} which shows   finite jet   determination. Nevertheless, we explain why, following his proof, we can not conclude $2$-jet determination in the general case. However we obtain $2$-jet determination (using the same notion of nondegeneracy) in two special cases: in $\Bbb C^4$ without the assumption of smoothness through the technics of analytic discs (see \cite {be-bl-me}), and for  codimension $2$  in $\Bbb C^N$ through the technics developed in \cite{KMZ} (see  \cite{Bl-Me}).

We point out that
finite jet determination problems for   submanifolds has attracted  much attention.
 We  refer  in particular  to the papers of Zaitsev \cite{Za}, Baouendi, Ebenfelt and Rothschild \cite{BER1}, Baouendi, Mir and Rothschild \cite{BMR}, Ebenfelt, Lamel and Zaitsev \cite{eb-la-za},  Lamel and Mir \cite{la-mi}, Juhlin \cite{ju}, Juhlin and Lamel \cite{ju-la}, Mir and Zaitsev \cite{mi-za} in the real analytic case, Ebenfelt \cite{eb}, Ebenfelt and Lamel \cite{eb-la}, Kim and Zaitsev \cite{ki-za}, Kolar, the second author and Zaitsev \cite{KMZ}  in the
$ \mathcal{C}^\infty$ case,  Bertrand and the first author \cite{be-bl}, Bertrand, the first and the second authors \cite {be-bl-me}  in the finitely smooth case.
\bigskip

The paper is organized as follows. In Section 2, we present the main protagonists and look at the model case. Section 3 is devoted to the various definitions of nondegeneracy. Finally, we give in Section 4 a detailed proof of the first part of the  main theorem of {\cite{Be2} and discuss the $2$-jet determination issue.

\section{Preliminaries}
\subsection{General issue} 
Let $M \subseteq \C^{N}$ be a $\mathcal{C}^{3}$ real submanifold of (real) codimension $d\ge 1$ through $p$. 

Some situations have to be excluded in order to get a finite dimensional stability group. For instance, the stability group of a complex submanifold of positive codimension is infinite dimensional: indeed, such a submanifold is given by $w=0$, and thus $(z,w)\mapsto (f(z),w)$ is a germ of automorphism at the origin as soon as $f(0)=0$ and $f$ is locally invertible at 0. At the other extreme, $\Aut(\R^N,0)$ is also infinite dimensional, because any mapping $F$ whose components are convergent power series with real coefficients, with $F(0)=0$ and the differential $dF_0$ is invertible, is a germ of automorphism. 

To avoid these situations, we deal with a submanifold $M$ that is generic, and thus CR, assuming the CR bundle is of positive dimension. Under these hypotheses, up to a linear transformation, then using the implicit function theorem, we can
assume that $p=0$ and the tangent space to $M\subseteq \C^n_z\times\C^d_w$ at 0 is $\Imm w_1=\cdots=\Imm w_d=0$,
hence $M$ is given by 
\begin{equation}\label{eq00}
\Imm w=\left(\begin{array}{c}{\,}_{\rm{real\ quadratic\ terms\ in}}\\  
    \Ree w,\Ree z,\Imm z\end{array}\right)+ O(|(\Ree w, z)|^3)
\end{equation}
Instead of grouping terms together according to their degree, we regroup them
by weight (see \cite{CM}): a weight 2 is assigned to $w$ and
a weight 1 is assigned to $z$. We get the $k$-th component of the second member of (\ref{eq00})
decomposed in three parts:
\begin{itemize}
\item terms of weight 2: a real quadratic form in $z$, that is
  $$q_k(z)+\overline{q_k(z)}+(Hermitian\ form\ in\ z)$$ with $q_k$ being
  a complex quadratic form;
\item terms of weight at least 3: linear combination of $(\Ree w_i)^2$, $(\Ree w_i)(\Ree z_j)$ and
  $(\Ree w_i)(\Imm z_j)$ for $1\le i\le d$, $1\le j\le n$;
\item $O(|(\Ree w, z)|^3)$, which would also consist of terms of weight at
  least 3 if the submanifold were real analytic.
\end{itemize}
After the local change of coordinates $z=z'$, $w_k=w_k'-2iq_k(z)$ ($1\le k\le d)$, the only remaining terms of weight 2 are Hermitian, and $M$ is given locally by the following system of equations:
\begin{equation}\label{equationreduite}
\begin{cases}
\Imm w_1= \transp\bar z A_1 z+ \left(\begin{array}{c}{\!}_{linear\ combination\ in}\\ {\ }_{(\Ree w_i)^2,\ z_j\Ree w_i,\ \overline{z}_j\Ree w_i}\end{array}\right)+O\ (|(z,\Ree w)|^3)\\
\ \ \ \ \vdots \\
\Imm w_d= \transp\bar z A_d z + \left(\begin{array}{c}{\!}_{linear\ combination\ in}\\ {\ }_{(\Ree w_i)^2,\ z_j\Ree w_i,\ \overline{z}_j\Ree w_i}\end{array}\right)+O\ (|(z,\Ree w)|^3)
\end{cases}
\end{equation}
where $A_1,\hdots,A_d$ are Hermitian matrices of size $n$.%, and the dots mean
%$$\textcolor{red}{\left(\begin{array}{c}{\!}_{linear\ combination}\\ {\ }_{of\ \Ree w_i\Ree w_j}\end{array}\right)+\left(\begin{array}{c}{\ }_{terms\ in\ z\ \Ree w\ }\\ {\ }_{and\ \bar{z}\ \Ree w}\end{array}\right)+o(|(z,\Ree w)|^3)}$$
\ Setting 
$$\langle z,z'\rangle=(\transp z A_1 z',\hdots,\transp z A_d z')$$
this system of equations can be written in a shortest way:
$$\Imm w=\langle \bar z,z\rangle + h(z,\Ree w)$$
where the Taylor expansion of $h$ at 0 consists of terms of weight at least 3.
We can also normalize the remaining term $h(z,\Ree w)$ by asking every term of its Taylor expansion at 0 to have no pure terms, that is to involve both $z$ and $\bar{z}$
(see \cite{BER} and \cite{Bo} section 7.2 for more details).

%\textcolor{blue}{En vertu du theoreme de Boggess, $h(z, \bar z, \Ree w)$ peut s'\' ecrire
%$$
%h(z, \bar z,\Ree w)= z^2 h_1(z, \bar z,\Ree w)+ {\bar z^2}h_2(z, \bar z,\Ree w)+ z\bar z h_3(z, \bar z,\Ree w)$$ $$ +\bar z (\Ree w )h_4(z, \bar z,\Ree w) +z(\Ree w) h_5(z, \bar z,\Ree w) + *({\Ree w}^2),
%$$ 
%avec $h_i(z, \bar z,\Ree w)\rightarrow 0$ if $(z, \bar z,\Ree w)\rightarrow 0.$ 
%Donc les termes sont soit des $o(|z|^2)$,
% soit des $o(\Ree w)$ 
 %ce qu'on \'ecrit $ o(|z|^2), \Ree w).$
%}

\subsection{The Levi map}
We recall that $A_1,\hdots,A_d$ are not uniquely determined but depend on the choice of the defining function of $M$. However, they contain some intrinsic geometric information on the submanifold because they are related to the Levi map, which is the generalization of the Levi form to the case of higher codimension. 

For $p\in M$, let $T_pM$ be the \emph{tangent bundle} to $M$ at $p$, $N_pM$ be the \emph{normal bundle} to $M$ at $p$ (that is, the orthogonal complement of $T_pM$ in $T_p(\R^{2N})$, and denote by $\pi_p:T_p(\R^{2N})\to N_pM$ the orthogonal projection. Following \cite{Bo}, we then define the Levi map for all $p\in M$:

\begin{df}
 The \emph{Levi map} of $M$ at $p$ is $L_p:T_p^{(1,0)}M\to N_pM$ defined on the space of holomorphic vectors tangent to $M$ at $p$ by 
$$L_p(X_p)=\frac{1}{2i}\pi_p(J[\overline{X},X](p))$$
where $X$ is any holomorphic vector field such that $X(p)=X_p$ and $J$ is the standard complex structure.
\end{df}

\

Note that the definition is valid for any CR submanifold of class $\mathcal{C}^2$, and that it is preserved under the action of CR-diffeomorphisms: let $M$, $M'$ be CR submanifolds and $L^M$, ${L^{M'}}$ their respective Levi maps, then if $F:M\to M'$ is a CR-diffeomorphism we get
\begin{equation}\label{LeviBiholo}
\forall p\in M,\ \forall X_p\in T_p^{1,0}M,\ L^M_p(X_p)=dF_p^{-1}({L^{M'}_{F(p)}}(dF_p\cdot X_p))
\end{equation}
Each component of $L_p$ being a Hermitian form, it is given by a unique sesquilinear form. We denote by $\mathcal{L}_p$ the corresponding sesquilinear map from $T_p^{(1,0)}M\times T_p^{(1,0)}M$ to $\C\otimes N_pM$, that is:
$\forall X_p,Y_p\in T_p^{(1,0)}M$,
\begin{eqnarray*}
2\mathcal{L}_p(X_p,Y_p)&=&\frac{1}{i}\pi_p(J[\overline{X},Y](p))\\
&=&(L_p(X_p+Y_p)-L_p(X_p)-L_p(Y_p))\\
 & &\ \ \ \ \ \ \ \ \ \ \ +i(L_p(X_p)+L_p(Y_p)-L_p(X_p+iY_p))
\end{eqnarray*}
where $\pi_p$ is extended between $\C\otimes T_p(\R^{2N})$ and $\C\otimes N_pM$.

\

In case $M$ is given by (\ref{equationreduite}), we obtain a convenient expression in coordinates for the Levi map of $M$ at 0 (\cite{Bo}, section 10.2):
\begin{equation}\label{Levimapscoord}
\forall z\in\Bbb C^n,\ L_0(z)=(\transp\bar z A_1 z,\hdots,\transp\bar z A_d z)=\langle\bar z,z\rangle
\end{equation}
where $N_0M$ is identified with $\R^d$. Similarly, we get:
$$ \forall z,z'\in\Bbb C^n,\ \mathcal{L}_0(z,z')=(\transp\bar z A_1 z',\hdots,\transp\bar z A_d z')=\langle\bar z,z'\rangle $$

\subsection{The case of a quadric}
To sort various conditions of non-degeneracy, we begin with the simplest case of a quadric submanifold $Q$, for which there is only a Hermitian part in $z$:
\begin{equation}\label{equationquadrique}
\begin{cases}
\Imm w_1= \transp\bar z A_1 z\\
\ \ \ \ \vdots \\
\Imm w_d= \transp\bar z A_d z
\end{cases}
\end{equation}
It is easy to see that we need to put at least two conditions on the $A_j$ if we expect $\Aut(Q,0)$ to be finite dimensional. This  will  give  a motivation for the nondegeneracy condition in the sense of Beloshapka.

\

\begin{lem}{}\label{lemmecond}
If $\Aut(Q,0)$ is finite dimensional, then
\smallskip
\begin{center}
\begin{tabular}{cl}
$\aaa$ & the Hermitian matrices $A_1$,...,$A_d$ are linearly independent (equivalently on $\R$ or $\C$)\\
\ \\
$\bbb$ & $\bigcap_{j=1}^d\mathrm{Ker}A_j=\{0\}$
\end{tabular}
\end{center}
\end{lem}

\ \\

\begin{rem}\label{remdim}\ \\
- Condition $\aaa$ cannot be satisfied for $d>n^2$, the dimension of the space of Hermitian matrices of size $n$.\\ 
- Condition $\bbb$ is equivalent to say that if $\langle z,z'\rangle=0$ for all $z'$, then we get $z=0$.
\end{rem}

\begin{pf}
Assume first that $\aaa$ is not satisfied, that is for instance $\exists \lambda_1,\hdots,\lambda_{d-1}\in\R$ such that $A_d=\sum_{i=1}^{d-1}\lambda_iA_i$. After the complex linear change of coordinates
$$(z_1,\hdots, z_n,w_1,\hdots,w_d)\mapsto(z_1,\hdots,z_n,w_1,\hdots,w_{d-1},w_d-\sum_{i=1}^{d-1}\lambda_iw_i)$$ 
equation (\ref{equationquadrique}) becomes 
\begin{equation*}
\begin{cases}
\Imm w_1= \transp\bar z A_1 z\\
\ \ \ \ \vdots \\
\Imm w_{d-1}= \transp\bar z A_{d-1} z\\
\Imm w_d=0
\end{cases}
\end{equation*}
The stability group at 0 thus contains the following mappings 
$$(z_1,\hdots,w_d)\mapsto(z_1,\hdots,w_{d-1},f(w_d))$$ 
for $f(w_d)=\sum_{k=1}^{+\infty}a_kw_d^k$ being any convergent power serie with $a_1\not=0$ and $a_k\in\R$ for all $k$.

Assume now that $\bbb$ is not satisfied: there exists some $z$ of norm 1 in $\bigcap_{j=1}^d\mathrm{Ker}A_j$. We can always assume $z=(1,0,\hdots,0)$ up to a unitary transformation, that is, the first column of every $A_j$ is zero, and so is the first raw. This means that $\langle z,z\rangle=(\transp\bar z A_1 z,\hdots,\transp\bar z A_d z)$ does not depend on $z_1$, hence the stability group at 0 contains mappings as
$$(z_1,\hdots,w_d)\mapsto(f(z_1),z_2,\hdots,z_n,w_1,\hdots,w_d)$$ 
where $f$ is any holomorphic function in a neighborhood of 0 such that $f(0)=0$ and $f'(0)\not=0$.
\end{pf}

\section{Non-degeneracy in the sense of Beloshapka}
Lemma \ref{lemmecond} leads to the following.

\subsection{The definition and its properties}
\begin{df} 
The generic $\mathcal{C}^{3}$ real submanifold $M$ of $\C^N$ given by (\ref{equationreduite}) is said \emph{Levi non-degenerate} at 0 (in the sense of Beloshapka) if conditions $\aaa$ and $\bbb$ are both satisfied. 
\end{df}

This convenient definition was introduced in \cite{Be1,Be2}. We will check later that, even if it seems to depend on the defining function of $M$, the Levi non-degeneracy is actually a biholomorphic invariant. We first investigate the relations between $\aaa$ and $\bbb$.

In the hypersurface case $d=1$, $\aaa$ means that the only matrix $A_1$ is non zero, and $\bbb$ means it is invertible; thus $\bbb$ is exactly the usual Levi non-degeneracy condition for hypersurfaces and it obviously implies $\aaa$. Note that $\aaa$ and $\bbb$ are equivalent only for a real hypersurface in $\C^2$ ($n=d=1$).

\begin{lem}
In codimension $d\ge 2$, 

- $\bbb$ never implies $\aaa$.

- $\aaa$ needs $d\le n^2$. If moreover $(n-1)^2<d$, then $\aaa$ implies $\bbb$; otherwise the implication is always false.
\end{lem}

\begin{pf}
  The first point is obvious, since one can choose $A_1$ invertible and $A_2=\cdots=A_d=0$. The beginning of the second point comes from Remark \ref{remdim}. If $\aaa$ is satisfied and $\bbb$ is not satisfied: then, up to a change of basis, the first raw and the first column of the $A_i$ are zero, and the dimension $d$ of the vector space generated by the $A_i$ is thus less than the dimension of Hermitian matrices of size $(n-1)$, that is, $(n-1)^2$.

  If $d\le(n-1)^2$: then $n\ge 2$, and we can choose $d$ linear independent Hermitian matrices $B_1,\hdots,B_d$ of size $n-1$. Set $A_i=\left(\begin{array}{cc} B_k&\begin{array}{c}0\\ \vdots\end{array} \\ \begin{array}{cc}0& \hdots\end{array}&0\end{array}\right)$: the $A_i$ are linearly independent but $\bbb$ is not satisfied.
\end{pf}
  
Thus, what are the first cases to explore for $\aaa$ and $\bbb$ being satisfied? In $\C^2$ resp. $\C^3$, it remains only the hypersurface case ($n=d=1$ or $n=2,d=1$); in $\C^4$, the hypersurface case ($n=3,d=1$) or $n=d=2$. So {\it the first interesting situation for higher codimension concerns the real codimension 2 in $\C^4$}. We observe that in this case we just have to check $\aaa$, since $\aaa$ implies $\bbb$ (because $d>(n-1)^2$).

\subsection{Other non-degeneracy conditions}
Among the wide literature concerning finite jet determination problems, various (non equivalent) generalizations in higher codimension of the notion of non-degeneracy appear:
\begin{itemize}
\item non-degeneracy in \cite{Za,BER} is exactly condition $\bbb$;
\item Levi non-degeneracy in the sense of Tumanov \cite{Tu} means, using (\ref{Levimapscoord}), that there exists a (real) linear combination $\sum \lambda_jA_j$ that is invertible; or equivalently, that the conormal bundle $N^*M$ is totally real at any point $(0,\sum\lambda_j\partial\rho_j(0))$ where $(\rho_j)$ is a set of defining functions for $M$; this obviously implies $\bbb$ (but not $\aaa$, except in the hypersurface case, since one can choose $A_1$ invertible and $A_2=\hdots=A_d=0$).
\end{itemize}

Note that the Levi non-degeneracy in the sense of Beloshapka does not give the Levi non-degeneracy in the sense of Tumanov:
 
\begin{exa}\cite{Be2}
In $\C^6$, consider the following codimension 3 real submanifold:
\begin{equation*}
\begin{cases}
\Imm w_1=|z_1|^2\\
\Imm w_2=2\Ree(z_1\bar{z}_2)\\
\Imm w_3=2\Ree(z_1\bar{z}_3)
\end{cases}
\end{equation*}
that is 
$$A_1=\left(\begin{array}{ccc}1&0&0\\0&0&0\\0&0&0\end{array}\right)\qquad 
A_2=\left(\begin{array}{ccc}0&1&0\\1&0&0\\0&0&0\end{array}\right)\qquad
A_3=\left(\begin{array}{ccc}0&0&1\\0&0&0\\1&0&0\end{array}\right)$$
Conditions $\aaa$ and $\bbb$ are satisfied, but no linear combination of $A_1,A_2,A_3$ is invertible.
\end{exa}

\subsection{Geometric characterizations}
Conditions $\aaa$ and $\bbb$ have been defined in terms of a reduced equation of $M$. But they are equivalent to some geometric conditions, proving actually they are biholomorphically invariant. More precisely:

\begin{prop}
The following statements are equivalent to condition $\aaa$:
\begin{enumerate}[i)]
\item the image of $L_p:T_p^{1,0}M\to N_pM$ is not included in a hyperplan;
\item the \emph{Levi cone of $M$ at $p$} (that is, the linear hull of the image of $L_p$ in $N_pM$) has a non empty interior;
\item $M$ is of finite type at $p$ with $2$ the only H\"ormander number.
\end{enumerate}
\end{prop}

\begin{rem} 
According to (\ref{LeviBiholo}), the first point implies $\aaa$ is biholomorphically invariant. The second point means 
$M$ has {\it a generating Levi form} at $p$ in the sense of \cite{Tu}.
\end{rem}

\begin{pf}
\begin{enumerate}[i)]
\item Expressing the Levi map in convenient coordinates (\ref{Levimapscoord}),
 $$L_p(T_p^{(1,0)}M)=\{(\transp\bar z A_1 z,\hdots,\transp\bar z A_d z)\ |\ z\in\C^n\}\subset\R^d$$
  is included in a hyperplane if and only if there exists $(\lambda_1,\hdots,\lambda_d)\in\R^d\setminus\{(0,\hdots,0)\}$ such that 
$$\forall z\in\C^n,\ \lambda_1\transp\bar z A_1 z+\hdots+\lambda_d\transp\bar z A_d z=0$$
  \textsl{i.e.} the Hermitian matrix $\lambda_1A_1+\hdots+\lambda_dA_d$ is zero.
\item Let $\Gamma_p$ be the Levi cone of $M$ at $p$: by definition, $\Gamma_p$ is included in a hyperplane if and only if $L_p(T_p^{(1,0)}M)$ is included in a hyperplane. But since $\Gamma_p$ is convex, it is included in a hyperplane if and only if its interior is empty.
\item Recall that $M$ being \emph{of finite type at $p$ with $2$ the only H\"ormander number} means that there exist CR vector fields $L_1,\hdots,L_j$ such that $L_1,\hdots,L_j,\overline{L_1},\hdots,\overline{L_j}$ and their Lie brackets generate the complexified tangent space $\C\otimes T_pM$. Or, equivalently, that there exists a basis of vector fields of $T^{(1,0)}M$: $X_1,\hdots,X_n$, such that $X_1(p),\hdots,X_n(p)$, $\overline{X_1(p)},\hdots,\overline{X_n(p)}$ and the $[\overline{X_i},X_j](p)$ generate $\C\otimes T_pM$.\\
By definition of the sesquilinear Levi map, this happens if and only if $Im(\mathcal{L}_p):=\{\mathcal{L}_p(X_p,Y_p)\ |\ X_p,Y_p\in T_p^{(1,0)}M\}$ generates $\C\otimes N_pM$. Since $\mathcal{L}_p(X_p,Y_p)$ can be written as a complex linear combination of values of $L_p$, we get
$Im(L_p)\subset Im(\mathcal{L}_p)\subset \mathrm{Span}_\C(Im(L_p))$ and thus $\mathrm{Span}_\C(Im(\mathcal{L}_p))=\mathrm{Span}_\C(Im (L_p))=\C\otimes\mathrm{Span}_\R(Im (L_p))$. Finally, $\mathrm{Span}_\C(Im(\mathcal{L}_p))=\C\otimes N_pM$ if and only if $\mathrm{Span}_\R(Im (L_p))=N_pM$.
\end{enumerate}
\end{pf}

\begin{rem}
  If the sesquilinear Levi map $\mathcal{L}_p$ of $M$ is surjective, then the first condition holds. But in contrary of what is said in Lemma 4.3.26 page 108  of \cite{BER}, the surjectivity of the sesquilinear Levi map $\mathcal{L}_p$ is not equivalent to $\aaa$, as we  will show in the following example.
\end{rem}

\begin{exa}
Take  in $\C^6,$
\begin{equation*}
\begin{cases}
\Imm w_1=|z_1|^2\\
\Imm w_2=|z_2|^2\\
\Imm w_3=2\Ree (z_1\bar z_2)\\
\Imm w_4=2\Imm (z_1 \bar z_2)\\
\end{cases}
\end{equation*}

Indeed, 
$$ A_1=
\begin{pmatrix}
1&0\\
\,\,0&0
\end{pmatrix},\quad
A_2=
\begin{pmatrix}
0&0\\
\,\,0&1
\end{pmatrix},
\quad
 A_3=
\begin{pmatrix}
0&1\\
\,\,1&0
\end{pmatrix},\quad
A_4=
\begin{pmatrix}
0&i\\
\,\,-i&0
\end{pmatrix}
$$
They are linearly independent but
$\mathcal{L}_0 : \C^2\times\C^2\to\C^4$ is expressed by
$$\mathcal{L}_0(z,z')=(\bar z_1 z'_1, \  \bar z_2 z'_2, \ \bar z_1 z'_2 + \bar z_2z'_1,\  i\bar z_1 z'_2 -i\bar z_2 z'_1)$$
and after a composition with a linear invertible map, we obtain
$$(z,z')\mapsto(\bar z_1 z'_1, \  \bar z_2 z'_2, \ \bar z_1 z'_2 ,\  \bar z_2 z'_1)$$
This map is not surjective since  it takes its values in 
$$\{ (t_1, t_2, t_3, t_4) \in \C^4 \ | \ t_1t_2= t_3t_4 \}$$
\end{exa}

\begin{rem}
Note that the surjectivity of $\mathcal{L}_p$ does not imply the surjectivity of $L_p$, even if $\aaa$ and $\bbb$ are satisfied: for instance, with $d=n=2$,
$$A_1=\left(\begin{array}{cc}1&0\\0&0\end{array}\right)\quad\mathrm{et}\quad A_2=\left(\begin{array}{cc}0&0\\0&1\end{array}\right)$$
    The image of $\mathcal{L}_0$ is $\C^2$ but the image of $L_0$ is $(\R^+)^2$.
\end{rem}

\

\begin{rem} Condition $\bbb$ is biholomorphically invariant since it means that 
$$\left(\forall X\in T_p^{(1,0)}M,\ \mathcal{L}_p(Y,X)=0\right)\quad \Longrightarrow\quad Y=0$$ 
\end{rem}
We recall the following definition.
\begin{df}\cite{BER}
A smooth generic submanifold $M \subset \Bbb C^N$ is holomorphically nondegenerate (in the formal sense) at $p \in M$ if there is no nontrivial formal holomorphic vector field  at $p$ tangent to $M.$
\end{df}
\begin{rem}
Note that if $M$ is real analytic and  holomorphically nondegenerate (in the formal sense) at $p,$ then $M$ is holomorphically nondegenerate at $p$ in the classical sense, that is, there is no nontrivial  holomorphic vector field  at $p$ tangent to $M.$  See Proposition 11.7.4 in \cite{BER}.
\end{rem}
 The following proposition shows that condition $\bbb$  has also some geometric meaning.
\begin{prop}
 Condition $\bbb$   implies that $M$  is holomorphically non degenerate (in the formal sense), provided that $M$ is smooth.
\end{prop}

\begin{pf}
Using  Proposition 2.1 in \cite{Za} and Theorem 11.5.1 in \cite{BER}, we deduce that $Q$ given by \eqref{equationquadrique} is holomorphically nondegenerate, and hence $M$ is  holomorphically nondegenerate. 
\end{pf}

\section{Discussing the Theorem of Beloshapka}
Let $M$ be a generic smooth real submanifold in $\C^N.$ The following statement obtained in $1989$ by Beloshapka shows that the  non-degeneracy condition  in the sense of Beloshapka was the good hypothesis for finite-jet determination  for  $M$ of finite type   with $2$ the only H\"ormander number. Later on, in 1997,  Zaitsev obtained the bound  $2(1+ \text{codim} \  M)$ for  the number of derivatives needed in the case where  $M$ is  real-analytic. (See Theorem 1.2 in \cite{Za}).

\begin{thm}\cite{Be2}\label{THM}\\
Let $M$ be a generic smooth real submanifold in $\C^N$, and $q\in M$. Assume $M$ is nondegenerate in the sense of Beloshapka at $q$: then there exists some $k$, only depending on the Levi map of $M$ at point $q$, such that the elements of $Aut(M,q)$ are uniquely determined by the values of the derivatives up to order $k$ at point $q$.
\end{thm}

Note that the result was stated for real analytic submanifolds, but the proof given in \cite{Be2} actually works for smooth ($\mathcal{C}^\infty$) submanifolds, but does not give a precise estimate of $k$. In this section, we fulfill the details of this original proof. As before, we denote by $d$ the real codimension of $M$, described by (\ref{equationreduite}):
$$v=\langle\bar z,z\rangle+h(z,u)$$
where $z\in\C^n$, $w=u+iv$ with $u,v\in\R^d$, and $\langle\bar z,z'\rangle=(\transp\bar z A_1 z',\hdots,\transp\bar z A_d z')$. Assume the Hermitian matrices $A_1,\hdots,A_d$ satisfy conditions $\aaa$ and $\bbb$. The proof of Theorem \ref{THM} relies on the following step.

\subsection{The basic identity}
Let $(F,G):\C^{n+d}\to\C^n\times\C^d$ be in $\mathrm{Aut}(M,0)$. Since $(F,G)$ maps $M$ into $M$, we get 
\begin{eqnarray*}
\lefteqn{\Imm G(z,u+i[\langle\bar z,z\rangle+h(z,u)])}\\
 &=&\left\langle\overline{F(z,u+i[\langle\bar z,z\rangle+h(z,u)])},F(z,u+i[\langle\bar z,z\rangle+h(z,u)])\right\rangle\\
 &\ &+\ h(F(z,u+i[\langle\bar z,z\rangle+h(z,u)]),\Ree G(z,u+i[\langle\bar z,z\rangle+h(z,u)]))
\end{eqnarray*}
We first follow \cite{CM}: this equality between two smooth functions leads to the equality of their Taylor expansions.
Let us write $F(z,w)=az+bw+\hdots$ and $G(z,w)=\alpha z + \beta w +\hdots$ where $a,b,\alpha,\beta$ are complex linear maps, and the dots represent higher order terms. Injecting in the previous equality, we identify the first coefficients:
\begin{itemize}
\item the right member contains no linear term, so $\alpha=0$ and $\Imm\beta=0$;
\item the Hermitian terms with respect to $z$ give $\beta\langle\bar{z},z\rangle=\langle\overline{az},az\rangle$.
\end{itemize}

Assigning weight one to $z$ and weight two to $w$, we decompose $F$ and $G$ into weighted homogeneous polynomials:
$$F=\sum_{q=0}^{+\infty}F_q\quad,\quad G=\sum_{q=0}^{+\infty}G_q$$
where  $F_q(tz,t^2w)=t^qF_q(z,w)$ and the same for $G_q$. We thus know that
$$\left\{\begin{array}{l} 
F_0=0\ ,\ F_1=az\\
G_0=0\ ,\ G_1=0\ ,\ G_2=\beta w\ \mathrm{where}\ \Imm\beta=0\ \mathrm{and}\ \beta\langle\bar{z},z\rangle=\langle\overline{az},az\rangle
\end{array}\right.$$
(notice that, since $(F,G)$ is biholomorphic at the origin, $a$ and $\beta$ must be invertible).

\bigskip

In the same way, for $q>2$, isolating the $q$-th component in the left member gives
$$\Imm G_q(z,u+i\langle\bar z,z\rangle) + \mathrm{terms\ in}\ G_r\ \mathrm{with}\ r<q$$
and in the right member
\begin{equation*}
\left\langle\overline{F_{q-1}(z,u+i\langle\bar z,z\rangle)},az\right\rangle+\left\langle\overline{az},F_{q-1}(z,u+i\langle\bar z,z\rangle)\right\rangle + \mathrm{terms\ in}\ F_{r-1}\ \mathrm{with}\ r<q 
\end{equation*}
So we obtain that the sequences $(F_q)$ and $(G_q)$ satisfy the following recurring system:

\begin{equation}\label{avec2dmb}
\forall q>2,\ \Ree\left(iG_q+2\langle\overline{F_{q-1}}, az\rangle\right)_{|v=\langle\bar z,z\rangle}=\mathrm{function\ of}\ G_r,F_{r-1}\ \mathrm{with}\ r<q
\end{equation}

\bigskip

Assume that $F_0,\hdots,F_{q-2}$ and $G_0,\hdots,G_{q-1}$ are known. Then the difference between two solutions $(F_{q-1},G_q)$ and $(\tilde{F}_{q-1},\tilde{G}_q)$ of (\ref{avec2dmb}) is a weighted homogeneous polynomial solution of the equation $\Ree(ig+2\langle \bar{f}, az\rangle)_{|v=\langle\bar{z},z\rangle}=0$. %We will prove that the only holomorphic solutions of this homogeneous equation are polynomials of degree less or equal to 2, and thus are of weight less or equal to 4. In particular, a solution of the recurring system (\ref{avec2dmb}) is uniquely determined by the initial values $F_0,\hdots,F_3$ and $G_0,\hdots,G_4$, and more precisely by the terms of degree less or equal to 2 of $F_0,\hdots,F_3$ and $G_0,\hdots,G_4$. This means that if two elements of $\Aut(M,0)$ have the same 2-jet at the origin, they are equal.
We will prove that the polynomial solutions of this homogeneous equation are of weight $q$ bounded by some $q_0$. In particular, a solution of the recurring system (\ref{avec2dmb}) is uniquely determined by the initial values $F_0,\hdots,F_{q_0-1}$ and $G_0,\hdots,G_{q_0}$. Let $k$ be the maximal degree appearing in terms of weight $q_0$: then if two elements of $\Aut(M,0)$ have the same $k$-jet at the origin, they are equal.

\medskip

Thus we aim to prove that (weighted homogeneous) polynomial solutions of the equation $\Ree(ig+2\langle \bar{f}, az\rangle)_{|v=\langle\bar{z},z\rangle}=0$ are of bounded degree. Using that $a$ and $\beta$ are invertible linear maps and that $\beta$ is real satisfying $\beta\langle\bar{z},z\rangle=\langle\overline{az},az\rangle$, we may replace $f$ by $af$ and $g$ by $\beta g$ and get
\begin{equation}\label{condequ}
\Ree(ig+2\langle\bar{f}, z\rangle)_{|v=\langle\bar{z},z\rangle}=0
\end{equation}

\noindent Theorem \ref{THM} will follow from 

\begin{thm}
%If the basic identity given by (\ref{condequ}) has holomorphic solutions $(f,g)$, those are polynomials of degree less or equal to 2.
The polynomial solutions of (\ref{condequ}) are of bounded degree. More precisely, the partial degree in the $z_j$ variables is bounded by 2.
\end{thm}  

\subsection{Towards a system of PDE}
Assume that $(f,g)$ is a holomorphic solution of (\ref{condequ}). We decompose $f=\sum_{p=0}^{+\infty}f_p$ and $g=\sum_{p=0}^{+\infty}g_p$ in the following way:
$$f_p(z,w)=\sum_{J}a^{(p)}_J(z)w^J\ ,\ g_p(z,w)=\sum_{J}b^{(p)}_J(z)w^J$$
where $a^{(p)}_J(z)$, $b^{(p)}_J(z)$ are holomorphic homogeneous polynomials in $z$ of degree $p$. We may then rewrite (\ref{condequ}) as

\begin{equation*}
(\ref{condequ}')\ \ \ \ \ \ \ \Ree\left(i\sum_pg_p(z,u+i\langle\bar{z},z\rangle)+2\left\langle \sum_p\overline{f_p(z,u+i\langle\bar{z},z\rangle)}, z\right\rangle\right)\equiv 0\ \ \ \ \ \ 
\end{equation*}

In (\ref{condequ}'), we sort the terms by bidegree $(k,l)$ in $(z,\bar{z})$, that is, grouping the $z^K\bar{z}^L$ such that $|K|=k$ et $|L|=l$. 

\begin{enumerate}[(a)]
\item\label{0-0} terms of bidegree $(0,0)$:
$\Ree\left(i\sum_Jb^{(0)}_J(z)u^J+0\right)$ \textsl{i.e.} the $b^{(0)}_J$ are real and
$$\Imm(g_0(u))\equiv 0$$
\item\label{1-0} terms of bidegree $(1,0)$: $(ig_1(z,u)+2\langle \overline{f_0(u)},z)\rangle)+0$ because the conjugate $\overline{(ig_1(z,u)+2\langle\overline{f_0(u)},z\rangle)}$ does not contain such terms; so 
$$ig_1(z,u)+2\langle\overline{f_0(u)},z\rangle\equiv 0$$
\item\label{k-0} terms of bidegree $(k,0)$ for $k\ge 2$: $\Ree(ig_k(z,u))$; since $g_k$ is holomorphic, $\Ree(ig_k(z,u))\equiv 0$ implies $g_k=0$, thus
$$\forall k\ge 2,\ g_k=0$$
\item\label{1-1} terms of bidegree $(1,1)$: they only involve $g_0$ and $f_1$. Note that
\begin{eqnarray*}
  g_0(z,u+i\langle\bar{z},z\rangle)&=&\sum_{J}b^{(0)}_J(u+i\langle\bar{z},z\rangle)^J\\
   &=&\sum_{J=(j_1,\hdots,j_s)}b^{(0)}_J\prod_{s=1}^d(u_s+i\langle\bar{z},z\rangle_s)^{j_s}
\end{eqnarray*}
and the $(1,1)$-part is 
  $\sum_{J=(j_1,\hdots,j_s)}b^{(0)}_J\sum_{l=1}^dj_lu_l^{j_l-1}i\langle\bar{z},z\rangle_l\prod_{s\not=l}u_s^{j_s}$.
  We thus introduce the operator $\Delta$, associating to any differentiable map $\varphi(z,u)$ the derivative of $\varphi$ with respect to $u$ at point $(z,u)$ in the direction $\langle\bar{z},z\rangle$, that is
  $$\Delta\varphi (z,u)=\partial_u\varphi_{(z,u)}\cdot\langle\bar{z},z\rangle$$
Note that $\displaystyle\Delta^k\varphi(z,u)=\partial^k_u\varphi_{(z,u)}\cdot(\langle\bar{z},z\rangle,\hdots,\langle\bar{z},z\rangle)$ 
when it is well defined.
This gives  
$$\Ree(-\Delta g_0(z,u)+2\langle\overline{f_1(z,u)},z\rangle)\equiv0$$
\item\label{2-1} terms of bidegree $(2,1)$: in $i\sum_pg_p(z,u+i\langle\bar{z},z\rangle)+2\left\langle \sum_p\overline{f_p(z,u+i\langle z,z\rangle)}, z\right\rangle$, they only come from the $g_1$ and the $f_0$ parts; in the conjugate, only from the $f_2$ part. More precisely, we get
$$-\Delta g_1(z,u)+2\langle\bar{z},f_2(z,u)\rangle-2i\langle\overline{\Delta f_0(z,u)},z\rangle\equiv 0$$
\item\label{k+1-1} terms of bidegree $(k+1,1)$ for $k\ge 2$:
$$-\Delta g_k(z,u)+2\langle\bar{z},f_{k+1}(z,u)\rangle\equiv 0$$
\item\label{2-2} terms of bidegree $(2,2)$: now we are interested in the terms of bidegree $(2,2)$ in $(u+i\langle\bar{z},z\rangle)^J$, and they are given by $\frac{(i)^2}{2!}\Delta^2$. Consequently we obtain 
  $\Ree(-\frac{i}{2}\Delta^2g_0(z,u)+2\langle\overline{i\Delta f_1(z,u)},z\rangle)$, and thus
  $$\Imm(\Delta^2g_0(z,u)+4\langle\overline{\Delta f_1(z,u)},z\rangle)\equiv 0$$
\item\label{3-2} terms of bidegree $(3,2)$:
$$-\frac{i}{2}\Delta^2g_1(z,u)+2i\langle\bar{z},\Delta f_2(z,u)\rangle-\langle\overline{\Delta^2f_0(z,u)},z\rangle)\equiv 0$$  
\item\label{3-3} terms of bidegree $(3,3)$: as previously, the terms of bidegree $(3,3)$ in $(u+i\langle\bar{z},z\rangle)^J$ are given by applying $\frac{(i)^3}{3!}\Delta^3$, thus
$$\Ree(\Delta^3g_0-6\langle\overline{\Delta^2f_1},z\rangle)\equiv 0$$
\end{enumerate}

\bigskip 

\noindent From (\ref{k-0}) and (\ref{k+1-1}), we obtain for $k\ge 3$: $\langle\bar{z},f_k(z,u)\rangle\equiv 0$, or equivalently (since $f_k$ is holomorphic in $z$) that for all $(z,u)$, $\forall z',\ \langle z',f_{k+1}(z,u)\rangle =0$. According to condition $\bbb$, we have $f_{k+1}(z,u)=0$ thus
$$\forall k\ge 3,\ f_k=0$$
and so the only possibly non-zero parts are $g_0$, $g_1$, $f_0$, $f_1$, $f_2$.

\

- Eliminate $g_1$ in (\ref{2-1}):\\
by (\ref{1-0}), $g_1(z,u)=2i\langle\overline{f_0(u)},z\rangle$, so we have $\Delta g_1(z,u)=2i\langle\overline{\Delta f_0(z,u)},z\rangle$, and then (\ref{2-1}) becomes
$$ \langle\bar{z},f_2(z,u)\rangle\equiv 2i\langle\overline{\Delta f_0(z,u)},z\rangle$$

- Eliminate $g_1$ in (\ref{3-2}):\\
since $\Delta g_1(z,u)=2i\langle\overline{\Delta f_0(z,u)},z\rangle$, we obtain $\langle\bar{z},\Delta f_2(z,u)\rangle\equiv 0$ and thus
$$\langle\overline{\Delta^2 f_0(z,u)},z\rangle\equiv 0$$

- Eliminate $\Imm g_0$ in (\ref{2-2}):\\
 $\Imm g_0(u)=0$ means the $b^{(0)}_J$ are real, thus $\Imm(\Delta^k g_0(z,u))\equiv 0$ for all $k$. Then (\ref{2-2}) becomes
$$\Imm(\langle\overline{\Delta f_1(z,u)},z\rangle)\equiv 0$$

- Eliminate $\Imm g_0$ in (\ref{3-3}):\\ 
applying $\Delta^2$ to (\ref{1-1}) gives $2\Ree(\Delta^2\langle\overline{f_1(z,u)},z\rangle)=\Ree(\Delta^3 g_0(z,u))$ and by replacing in (\ref{3-3}), we get
$$\Ree(\Delta^3g_0(z,u))\equiv 0$$

\

Finally, we obtain that $(f,g)$ is a polynomial solution of (\ref{condequ}) if and only if the $f_k$ and $g_k$ form a polynomial solution of the following system:
\begin{equation}\label{system}
\begin{cases}
(8a)\ \ \! \forall k\ge 3,\ f_k=0\ \mathrm{and}\ g_{k-1}=0\\
(8b)\ \begin{cases}
ig_1+2\langle\bar{f_0},z\rangle=0\\
\langle\bar{z},f_2\rangle=2i\langle\overline{\Delta f_0},z\rangle\\
\langle\overline{\Delta^2f_0,}z\rangle=0
\end{cases}\\
(8c)\ \begin{cases}
\Imm g_0=0\\
2\Ree(\langle\bar{f_1},z\rangle)-\Ree(\Delta g_0)=0\\
\Imm\langle\overline{\Delta f_1},z\rangle=0\\
\Ree\Delta^3g_0=0
\end{cases}
\end{cases}
\end{equation}
where all maps are evaluated on $\{\Imm w=0\}$, that is at $(z,u)$.

\medskip

Note that $(\ref{system}a)$ exactly means that {\it the partial degree of $(f,g)$ in the $z_j$ variables is at most 2.}

\subsection{The solution space of (\ref{system}) is finite dimensional}
The maps $f_0$ and $g_0$ do not depend on $z$. We decompose $f_1$, $f_2$ and $g_1$ in the following way:
$$f_1(z,u)=\sum_{J}a^{(1)}(z)_Ju^J=\sum_J\left(\sum_{i=1}^n\alpha_{J,i}z_i\right) u^J=\sum_{i=1}^n\left(\sum_J\alpha_{J,i}u^J\right)z_i$$
thus we set $\displaystyle f_1(z,u)=\sum_{i=1}^n\phi_i(u)z_i$
and similarly $\displaystyle g_1(z,u)=\sum_{i=1}^n\psi_i(u)z_i$ and $\displaystyle f_2(z,u)=\sum_{i\le j}\Phi_{i,j}(u)z_iz_j$. Note that $\Delta f_1(z,u)=\sum_{i=1}^n\Delta\phi_i(z,u)z_i$, and the same for $g_1$ and $f_2$.

All the equalities in (\ref{system}b)-(\ref{system}c) are polynomials in $z$ and $\bar{z}$. If all their coefficients are equated to zero, then we get an overdeterminated homogeneous system of linear partial differential equations in the $n$ components of $\overline{f_0}$, $\Ree\phi_i$, $\Imm\phi_i$, $\Phi_{i,j}$ and the $d$ components of $\Ree g_0$, $\Imm g_0$, $\psi_i$ ($i,j=1,\hdots,n$) with {\it constant} coefficients:
\begin{equation}\label{EDP}
P(D)\Theta=0 \ \mathrm{on}\ \Omega\subset\R^d
\end{equation}
where the $m$ components of $\Theta$ are precisely the unknown functions, $D=(\frac{\partial}{\partial u_1},\hdots,\frac{\partial}{\partial u_d})$ , $P=(P_{r,s})$ is a $l\times m$-rectangular matrix with polynomial entries $P_{r,s}\in\C[X_1,\hdots,X_d]$, and $\Omega$ is an open convex neighborhood of the origin.

\medskip

As for linear ordinary differential equations with constant coefficients, we begin with looking at the exponential solutions. Set $\Theta(u)=e^{\zeta\cdot u}\cdot\theta$, where the vector $\theta$ is independent of $u$ and $\zeta\cdot u=\sum_{j=1}^d\zeta_ju_j$: it is a non-trivial solution of (\ref{EDP}) if and only if $P(\zeta)\cdot\theta=0$.
Let $V_P=\{\zeta=(\zeta_1,\hdots,\zeta_d)\in\C^d/\ \mathrm{rank}(P(\zeta))<m\}$ be the characteristic set of the system (\ref{EDP}). Note that $\zeta\in V_P$ exactly means that $\zeta$ induces a {\it complex valued} non-trivial exponential solution of (\ref{EDP}). 

\begin{lem}
The only possible characteristic value for the system (\ref{EDP}) is 0.
\end{lem}

\begin{pf}
Note that if we set $g_0^r(u)=\Ree (g_0)(u)$, $g_0^i(u)=\Imm (g_0)(u)$, $f_1^r(z,u)=\Ree (f_1)(z,u)=\sum_{i=1}^n\left(\sum_J\Ree(\alpha_{J,i})u^J\right)z_i$, $f_1^i(z,u)=\Imm (f_1)(z,u)$, then
$$(\ref{system}c)\Longleftrightarrow \begin{cases}
g_0^i=0\\
2\langle f_1^r(\bar{z},u)-if_1^i(\bar{z},u),z\rangle+2\langle\bar{z},f_1^r(z,u)+if_1^i(z,u)\rangle-\Delta(g_0^r)(z,u)=0\\
\langle\Delta(f_1^r)(\bar{z},u)-i\Delta(f_1^i)(\bar{z},u),z\rangle-\langle\bar{z},\Delta(f_1^r)(z,u)+i\Delta(f_1^i)(z,u)\rangle=0\\
\Delta^3g_0^r=0
\end{cases}$$ 
Assume that $\Theta(u)=e^{\zeta\cdot u}\cdot\theta$ is solution of (\ref{EDP}) with $\zeta\in\C^d$ and $\theta$ a complex vector. By definition of $\Theta$, this means we can write
$$\overline{f_0}(u)=e^{\zeta\cdot u}\cdot\tilde{f}_0\ ,\ f_2(z,u)=e^{\zeta\cdot u}\cdot\tilde{f}_2(z)\ , g_1(z,u)=e^{\zeta\cdot u}\cdot\tilde{g}_1(z)$$
$$g_0^r(u)=e^{\zeta\cdot u}\cdot\tilde{g}_0^r,\ g_0^i(u)=0,\ f_1^r(z,u)=e^{\zeta\cdot u}\cdot\tilde{f}_1^r(z),\ f_1^i(z,u)=e^{\zeta\cdot u}\cdot\tilde{f}_1^i(z)$$
Assume $\zeta\not=0$ and check that the non-degeneracy hypothesis implies $\tilde{f}_0=\tilde{f}_1^r=\tilde{f}_1^i=\tilde{f}_2=0$ and $\tilde{g}_0^r=\tilde{g}_0^i=\tilde{g}_1=0$. 

\medskip

Since $\overline{\Delta f_0(u)}=\Delta\overline{f_0}(u)=(\zeta\cdot \langle\bar{z},z\rangle)e^{\zeta\cdot u}\cdot\tilde{f}_0$, (\ref{system}b) becomes
$$\begin{cases}
i\tilde{g}_1(z)+2\langle\tilde{f}_0,z\rangle=0\\
\langle\bar{z},\tilde{f}_2(z)=0\\
(\zeta\cdot \langle\bar{z},z\rangle)^2\langle\tilde{f}_0,z\rangle=0
\end{cases}$$
Condition $\aaa$ gives $\zeta\cdot \langle\bar{z},z\rangle$ is not identically zero since $\zeta\not=0$, hence $\langle\tilde{f}_0,z\rangle=0$ and condition $\bbb$ implies $\tilde{f}_0=0$; thus $\tilde{g}_1=0$. Since $\tilde{f}_2$ is holomorphic, we also get by condition $\bbb$ that $\tilde{f}_2=0$.

\medskip

Similarly, for (\ref{system}c) we obtain $\tilde{g}_0^r=\tilde{g}_0^i=0$, and the two remaining equations become
$$\langle\tilde{f}_1^r(\bar{z})-i\tilde{f}_1^i(\bar{z}),z\rangle\pm\langle\bar{z},\tilde{f}_1^r(z)+i\tilde{f}_1^i(z)\rangle=0$$
Since $\tilde{f}_1^r$ and $\tilde{f}_1^i$ are holomorphic, by condition $\bbb$ we also get $\tilde{f}_1^r\pm i\tilde{f}_1^i=0$.

\end{pf}

We conclude by means of the following statement:

\begin{thm}\cite{Ob}
In the situation described above, the following assertions are equivalent:
\begin{enumerate}[a)]
\item $\C[X_1,\hdots,X_d]^{m}/P^T\cdot\C[X_1,\hdots,X_d]^{l}$ is finite-dimensional;
\item the solution space $S=\{y\in\mathcal{C}^\infty(\R^d,\C)\ |\ P(D)y=0\}$ is finite-dimensional;
\item the characteristic set $V_P$ is of dimension zero.
\end{enumerate}
\end{thm}

\begin{rem}
If these equivalent conditions are satisfied, then the solution space $S$ has a direct sum decomposition (integral representation):
$$S=\bigoplus_{\zeta\in V_P}\mathcal{P}(\zeta)e^{\zeta\cdot u}$$
where $\mathcal{P}(\zeta)\subset\C[X_1,\hdots,X_d]^{m}$ are finitely dimensional spaces of polynomials, depending on $P$. In the general situation, only assuming that $V_P\not=\C^d$, the Fundamental Principle of Ehrenpreis-Palamodov  states that every smooth solution is still the limit of a sequence of exponential-polynomial solutions (\cite{Ho}, Theorem 7.6.14).
\end{rem}

\medskip

Thus the solution space of (\ref{system}) is finite dimensional. In particular the degrees of the polynomial solutions are bounded by a single constant, depending only on the system, that is, on the Levi map of the submanifold at the origin. 

\subsection{Identifying the problem}
 Following the proof of Beloshapka, page 242,   the degree with respect to $u$  of the polynomial solution of  (\ref{system}b), $f_0$  (respectively $f_2$), is claimed to be at most $1$ (respectively $0$).
One  writes
\begin{equation}\label{dim1}
f_0= a_0 + a_1 + \dots,
\end{equation}
\begin{equation}\label{dim2}
f_2= A_0 +A_1 + \dots,
\end{equation}
each of $a_p$ (respectively $A_q$) beeing of degree $p$ (respectively $q$) with respect to $u.$ One writes
\begin{equation} \label{dim3}
a_p=\sum  a_{j_1}\dots _{j_p} u_{j_1}\dots u_{j_p},
\end{equation}
where the coefficients $a_{j_1}\dots _{j_p}$ are symmetric with respect to the indices.
\begin{equation} \label{dim4}
A_q=\sum  A_{j_1}\dots _{j_q} u_{j_1}\dots u_{j_q},
\end{equation}
where the coefficients $A_{j_1}\dots _{j_q}$ are symmetric with respect to the indices.
Rewriting  then  (\ref{system}b) with the help of \eqref{dim1},\eqref{dim2}, \eqref{dim3} and \eqref{dim4}, one obtains 
\begin{equation}\label{dim5}
<\bar z,A_{j_1}\dots _{j_{p-1}}>=2ip\sum_q<\bar z,z>_q <\overline{a_{j_1}\dots _{j_{p-1}q}},z>.
\end{equation}
\begin{equation}\label{dim6}
p(p-1)\sum_q<\bar z,z>_{q_{1}} <\bar z,z>_{q_2}
<\overline{a_{j_1}\dots _{{j_{p-1}}{q_1}{q_2}}},z>=0.
\end{equation}
Assume now that \eqref{dim5} and \eqref{dim6} have a nonzero solution for some $p>1,$ that is, there exists an homogeneous solution of (\ref{system}b).  Then Beloshapka proposes to set, for $s>0,$ 
\begin{equation}\label{dim7}
A_{l_1\dots l_sj_1}\dots _{j_{p-1}}=A_{j_1}\dots _{j_{p-1}}, \ \ a_{l_1\dots l_sj_1}\dots _{j_{p}}=a_{j_1}\dots _{j_{p}}
\end{equation}
Beloshapka   then claims  that the homogeneous  polynomial obtained by \eqref{dim7} gives a nonzero  homogeneous solution of (\ref{system}b) of arbitrarily degree, which would lead to a contradiction. The problem is that these  coefficients $A_{l_1\dots l_sj_1}\dots _{j_{p-1}}$ and 
$a_{l_1\dots l_sj_1}\dots _{j_{p}}$ are NOT symmetric with respect to the indices.
\begin{rem}
To our knowledge, the bound   $k= 2(1+ \text{codim} \  M)$ obtained by Zaitvev in \cite{Za} seems to be the best up to now for generic submanifolds that are real-analytic and  nondegenerate in the sense of Beloshapka.
\end{rem}

\end{document}